\documentclass[12pt]{amsart}
\usepackage{amscd,amssymb}
\usepackage[graph,frame,poly,arc]{xy}  
\usepackage[plainpages,backref,urlcolor=blue]{hyperref}

\theoremstyle{plain}
\newtheorem{thm}[subsection]{Theorem}

\newtheorem{prop}[subsection]{Proposition}
\newtheorem{cor}[subsection]{Corollary}

\theoremstyle{definition}
\newtheorem{rk}[subsection]{Remark}
\newtheorem{definition}[subsection]{Definition}
\newtheorem{ex}[subsection]{Example}

\numberwithin{equation}{section}
\newcommand{\OO}{{\mathcal O}}

\newcommand{\A}{{\mathcal A}}
\newcommand{\wJ}{\widehat{{J}}}

\newcommand{\E}{{\mathcal E}}

\newcommand{\LL}{{\mathcal L}}
\newcommand{\I}{{\mathcal I}}

\newcommand{\al}{{\alpha}}
\newcommand{\be}{{\beta}}

\newcommand{\C}{\mathbb{C}}

\newcommand{\PP}{\mathbb{P}}

\newcommand{\N}{\mathbb{N}}

\DeclareMathOperator{\Ext}{Ext}

\DeclareMathOperator{\Hom}{Hom}
\DeclareMathOperator{\rank}{rank}
\DeclareMathOperator{\im}{im}

\DeclareMathOperator{\codim}{codim}

\begin{document}

\title [Jumping lines of bundles of logarithmic vector fields]
{On the jumping lines  of  bundles of logarithmic vector fields along plane curves}

\author[Alexandru Dimca]{Alexandru Dimca}
\address{Universit\'e C\^ ote d'Azur, CNRS, LJAD and INRIA, France }
\email{dimca@unice.fr}

\author[Gabriel Sticlaru]{Gabriel Sticlaru}
\address{Faculty of Mathematics and Informatics,
Ovidius University
Bd. Mamaia 124, 900527 Constanta,
Romania}
\email{gabrielsticlaru@yahoo.com }



\begin{abstract} For a reduced curve $C:f=0$ in the complex projective plane $\PP^2$, we study the set of jumping lines for the rank two vector bundle $T\langle C \rangle $ on $\PP^2$,
whose sections are the logarithmic vector fields along $C$. We point out the relations of these jumping lines with the Lefschetz type properties of the Jacobian module of $f$ and with the Bourbaki ideal of the module of Jacobian syzygies of $f$. In particular, when the vector bundle $T\langle C \rangle $ is unstable, a line is a  jumping line if and only if it meets the 0-dimensional subscheme defined by this Bourbaki ideal, a result going back to Schwarzenberger.
Other classical general results by Barth, Hartshorne and Hulek  resurface in the study of this special class of  rank two vector bundles. 

\end{abstract}

\maketitle

\section{Introduction} \label{s0}

Let $C:f=0$ be a reduced curve of degree $d$ in $X=\PP^2$, $S=\C[x,y,z]$ be the polynomial ring with the usual grading, and $AR(f)$ be the graded $S$-module of Jacobian syzygies of $f$, see equation \eqref{eqAR} below. 
Let $E_C$ be the locally free sheaf on $X$ corresponding to the graded $S$-module $AR(f)$, and recall that
\begin{equation} \label{equa1} 
E_C=T\langle C \rangle (-1),
\end{equation}
where $T\langle C \rangle $ is the sheaf of logarithmic vector fields along $C$ as considered for instance in \cite{AD, DS14,MaVa}.
For a line $L$ in $X$, the pair of integers $( d_1^L, d_2^L)$ such that  $ d_1^L \leq d_2^L$ (resp. without the condition $ d_1^L \leq d_2^L$), and such that $ E_C|_L \simeq \OO_L(-d_1^L) \oplus \OO_L(-d_2^L)$ is called the (ordered) splitting type (resp. the unordered splitting type) of $E_C$ along $L$, see for instance \cite{FV1,OSS}. Unless we say the opposite, in this paper we use the (ordered) splitting type.
For a generic line $L_0$, the corresponding splitting type $( d_1^{L_0}, d_2^{L_0})$ is known to be constant, see Theorem \ref{thm1} (2) below.
A line $L$ in $X$ is called a jumping line of order $o(L)$ for $E_C$ or, equivalently, for $T\langle C \rangle $, if 
$$o(L):=d_1^{L^0}-d_1^{L}>0,$$
see for instance \cite[Section 5]{Ha80}.

When the graded $S$-module $AR(f)$ is free (equivalently, when $E_C$ splits as a direct sum of two line bundles on $X$),
which can be considered as the simplest case, then the corresponding curve is called {\it free}, a notion going back to K. Saito \cite{S}.
When the minimal resolution of the graded $S$-module $AR(f)$ is slightly more complicated, we get the {\it nearly free} curves considered in \cite{AD,B+,DStRIMS, MaVa}, see Definition \ref{defFC} below for details. 
For a free curve $C$, there are no jumping lines for $E_C$, while for a nearly free curve $C$ the jumping lines for $E_C$, if they exist, are of order 1 and form a line $\LL$ in the dual projective space $\PP(S_1)$, dual to the jumping point $P(C)\in X$ associated to $C$ by S. Marchesi and J. Vall\` es in \cite{MaVa}. In this note we study the set of jumping lines for $E_C$ for {\it any reduced plane curve} $C$.
As one {\it motivation} for this study, note that when $C$ is a line arrangement $\A$ in $X$, the question whether the combinatorics of $\A$ determines the generic splitting type of the corresponding vector bundle $T\langle \A \rangle $ is actively considered, see for instance \cite[Question 7.12]{CHMN} or the relation to Terao's Conjecture on the freeness of line arrangement explained in \cite{AD}, see also Remark \ref{rktop}.

%\medskip

In section 2 we start by recalling  some basic notions and results, in particular Theorem \ref{thm1} which determines completely the generic splitting type $( d_1^{L_0}, d_2^{L_0})$ in terms of the minimal degree $r=mdr(f)$ of a Jacobian syzygy and the degree $d$ of the curve $C$. The invariant $r$ also decides whether the vector bundle $E_C$ is stable: the stability holds if and only if $2r \geq d$, see \cite{Se} or the discussion below in section 2.

In section 3 we study the {\it Hilbert function} $\{k \mapsto h^1(\PP^2,E_C(k))\}$.
 Theorem \ref{thmB2}, which treats the case $E_C$ stable, is similar to, and can be obtained from, a result by Hartshorne, namely \cite[Theorem 7.4]{Ha78}.
The other main result, Theorem \ref{thmB1}, which shows that when $E_C$ is unstable, the behavior of the above Hilbert function is very different, is in our opinion completely new.
In particular, this gives a {\it partial strong Lefschetz property of the Jacobian module} $N(f)$ in the case $E_C$ unstable, see Corollary \ref{corLef}.

In section 4 we relate the integer $d_1^L$ to some Lefschetz type properties of the multiplication by an equation $\al_L$ of the line $L$, acting on the Jacobian module $N(f)$,
see Proposition \ref{propthmA}. Then we define and establish the first properties of the $k$-th {\it jumping locus}  $V_k(C)$ of the curve $C$, which consists of all lines $L$ in $X$ such that $d_1^L \leq k$, see Theorem \ref{thmB}. The main results in this section are Corollaries \ref{corB} and \ref{corB1}.
We note that the claim in Corollary \ref{corB1} (1) fits perfectly well with a general result of W. Barth about the pure 1-dimensionality of the set of jumping lines of stable rank 2 bundles with even first Chern class, see Remark \ref{rkB}. For stable rank 2 bundles with odd first Chern class, to get a similar result, K. Hulek has introduced the notion of jumping line of second kind,
and his results in \cite{Hu} lead again to a {\it partial strong Lefschetz property of the Jacobian module} $N(f)$ in the case $E_C$ stable, see Remark \ref{rkB}.

In section 5 we introduce the main new technical tool, namely the Bourbaki ideal $B(C,\rho_1) \subset S$ associated to the curve $C$ and to a minimal degree Jacobian syzygy $\rho_1$ for $f$, see Theorem \ref{thmC}. This allows us to present the vector bundle $E_C$ as an extension of the ideal of a codimension 2 locally complete intersection by a line bundle. The general construction of this type goes back to Serre \cite{Serre}, and it was widely used to construct rank 2 vector bundles on $\PP^n$ for $n \geq 3$, see \cite[Chapter 1, section 5]{OSS} and the many references given there. {\it The new point in our approach is the very explicit description of the ideal
$B(C,\rho_1) $}. When $r \leq d/2$ (resp. $r <d/2$), a line $L$ is not a jumping line if (resp. if and only if) it avoids the support of the subscheme $Z(C,\rho_1)$ of $\PP^2$ defined by the ideal $B(C,\rho_1) $, see Theorem \ref{thmD}, Theorem \ref{thmD2} and Corollary \ref{corD}. This result generalizes the result of S. Marchesi and J. Vall\` es in \cite{MaVa} concerning nearly free curves, and in fact goes back to  Proposition 10 in Schwarzenberger's paper \cite{Sch1}. When $r >d/2$, a line $L$ avoiding the support of the subscheme $Z(C,\rho_1)$ may be a jumping line, but it satisfies $d_1^L \geq d-r-1$, see Theorem \ref{thmBD}, and this lower bound seems to be strict in many cases.
The strong dependence of the ideal $B(C,\rho_1)$ and of the scheme $Z(C,\rho_1)$ of the choice of the syzygy $\rho_1$ is illustrated in Example \ref{ex2}. This is not a surprize, since
$Z(C,\rho_1)$ is exactly the zero locus subscheme of the section of the vector bundle $E_C(r)$
associated to the syzygy $\rho_1$, as explained in Remark \ref{rkC}.

We conclude with five examples in section 6. In the first one, we discuss the case of smooth curves, and point out in particular, that for $d=2d'+1$ odd, the geometry of the 
 jumping locus curve $V_{d'-1}(C)$ is quite interesting. This is a special case of 
Barth's result in 
\cite[Theorem 2 and Section 7]{Ba}, saying that a rank 2 stable vector bundle on $\PP^2$ with even second Chern class is determined by the associated net of quadrics, having the curve $V_{d'-1}(C)$ as its discriminant.

The other four examples discuss singular curves, satisfy all $d_1^{L_0}=2$ and hence $V_2(C)=\PP(S_1)$, the set of all lines in $\PP^2$. A quintic $C$ such that $E_C$ is semistable is considered in Example \ref{exS1}.
 In Example \ref{ex1}, $C$ is again a singular quintic, the first jumping locus $V_1(C)$ is a smooth conic, hence the 1-dimensional irreducible components of the jumping loci are not necessarily lines, and this is related to another general result by W. Barth, on the smoothness of some sets of jumping lines, see Remark \ref{rkE}.
In Example \ref{ex0}, $C$ is a Zariski sextic with 6 cusps on a conic, the first jumping locus $V_1(C)$ is the union of a line $\LL$ and two points, and hence it is not pure dimensional.
In Example \ref{ex2}, $C$ is another singular sextic, the first jumping locus $V_1(C)$ consists of 11 points, and  the 0-th
jumping locus $V_0(C)$ is one of the points in $V_1(C)$.  
In the last three examples, the corresponding vector bundles $E_C$ are stable, and hence the structure of the jumping loci can be rather subtle even in the class of stable rank two vector bundles of type $E_C$.

%\medskip

As shown by these examples, the jumping loci $V_k(C)$ for any plane curve $C$ can be determined explicitly using a Computer Algebra software, in our case we have used the package SINGULAR, see \cite{Sing}.

%\medskip

The first author thanks AROMATH team at INRIA Sophia-Antipolis for excellent working conditions, and  A. Beauville, L. Bus\'e, B. Mourrain and C. Pauly for stimulating discussions. We also thank J. Vall\`es, who has helped us with the proof of Theorem \ref{thmD2}.

\section{Preliminaries}

  For the coordinate ring 
$S=\C[x,y,z]$ and a graded $S$-module $M$, let $M_k$ denote the homogeneous 
degree $k$-part of $M$ and, for an integer $m$, define the shifted graded $S$-module $M(m)$ by the condition $M(m)_k=M_{m+k}$ for any $k$.
 For $g \in S$, let $g_x,g_y,g_z$ denote the partial 
derivative of $g$ with respect to $x,y,z$. Then the graded $S$-module $AR(f)=AR(C) 
\subset S^{3}$ of 
\textit{all Jacobian relations} for $f$ is defined by 
\begin{equation} \label{eqAR} 
AR(f)_k:=\{(a,b,c) \in S^{3}_k \mid a f_x+b f_y+c f_z=0\}.
\end{equation}
Its sheafification $E_C:=
\widetilde{AR(f)}$ is a rank two vector bundle on $\PP^2$, see 
\cite{AD,S, Se} for details. More precisely, one has $E_C=T\langle C \rangle (-1)$,
where $T\langle C \rangle $ is the sheaf of logarithmic vector fields along $C$ as considered for instance in \cite{AD, DS14,MaVa}. We set 
$$ar(f)_m= \dim AR(f)_m=\dim H^0(\PP^2,E_C(m))$$
 for any integer $m$.
We have the following, see \cite{AD, DStRIMS}.

\begin{definition}
\label{defFC}

(1)\,\,
A curve $C$ is \textit{free} if the graded $S$-module
$AR(f)$ is free, say with  a  basis $\rho_1,\rho_2$.  If
$\deg \rho_i=d_i\ (i=1,2)$,
the multiset of integers $(d_1,d_2)$ is called the 
\textit{exponents} of a free curve $C$.

(2)\,\,
A curve $C$ is \textit{nearly free} 
 if the graded $S$-module
$AR(f)$ has a minimal generator system of syzygies $\rho_1, \rho_2,\rho_3$, such that 
the degrees $\deg \rho_i$ satisfy
$d_1 \le d_2=d_3$ and there is
 a relation 
$$
h \rho_ 1+\ell_2 \rho_2+\ell_3 \rho_3=0,
$$
for $h \in S$ and independent linear forms $\ell_2,\ell_3\in S$. The multiset $(d_1,d_2)$ is called the
\textit{exponents} of a nearly free curve $C$.

\end{definition}

 Let  $mdr(f):=\min \{k \mid 
AR(f)_k \neq (0)\}$ be the minimal degree of a Jacobian syzygy for $f$. In this paper we assume that 
$
mdr(f) \ge 1,
$
unless otherwise specified. 
Let $N(f)= \wJ_f/J_f$, with $J_f$ the Jacobian ideal of $f$ in $S$, spanned by the partial derivatives $f_x,f_y$, $f_z$ of $f$, and $\wJ_f$ the saturation of the ideal $J_f$ with respect to the maximal ideal 
${\bf m}=(x,y,z)$ in $S$. The quotient module $N(f)$ coincides with $H^0_{\bf m}(S/J_f)$ and is called the {\it Jacobian module} of $f$, or of the plane curve $C$, see \cite{Se}.
The quotient $M(f)=S/J_f$ is called the {\it Jacobian algebra} of $f$, and we denote
$$m(f)_k = \dim M(f)_k$$
for any integer $k$.
Let $\nu(C)= \dim N(f)_{\lfloor T/2 \rfloor }$, where $T=3(d-2)$ and note that we have the following, see \cite[Theorem 1.2]{Drcc}.
\begin{thm}
\label{thmNU}
Let $C:f=0$ be a reduced plane curve of degree $d$ and let $r=mdr(f)$.
Then the following hold.
\begin{enumerate}
\item If $r < d/2$, then
$$\nu(C)=\tau_{max}(d,r)-\tau(C),$$
where $\tau_{max}(d,r)=(d-1)^2-r(d-1-r)$.
\item If $r \geq (d-2)/2$, then 
$$\nu(C)= \left \lceil   \  \frac{3}{4}(d-1)^2 \    \right \rceil  -\tau (C).$$
\end{enumerate}
\end{thm}
Here, for any real number $u$,  $\lceil     u    \rceil $ denotes the round up of $u$, namely the smallest integer $U$ such that $U \geq u$.
It is known that the curve $C:f=0$ is free (resp. nearly free) if and only if $\nu(C)=0$ (resp. $\nu(C)=1$), see  \cite{Dmax, DPop, DStRIMS}.
Recall the definition of the global Tjurina number 
$$\tau(C)= \sum _{p \in C} \tau(C,p)$$ 
of the curve $C$, where $\tau(C,p)$ is the Tjurina number of the singularity $(C,p)$, and the fact that $\tau(C)$ is the degree of the Jacobian ideal $J_f$. 
This numerical invariant of $C$ occurs in the following formulas for the Chern numbers of the vector bundle $T\langle C \rangle (k)=E_C(k+1)$, namely one has
\begin{equation}
\label{Chern}
c_1(T\langle C \rangle (k))=3-d+2k, \  c_2(T\langle C \rangle (k))=d^2-(k+3)d+k^2+3k+3-\tau(C),
\end{equation}
see for instance in \cite[Equation (3.2)]{DS14}. Associated to the vector bundle $E_C$ there is the {\it normalized} vector bundle $\E_C$, which is the twist of $E_C$ such that $c_1(\E_C) \in \{-1,0\}$. More precisely,
when $d=2d'+1$ is odd, then
\begin{equation}
\label{Chern1}
\E_C=E_C(d') \text{ and } c_1(\E_C)=0, \ c_2(\E_C)=3(d')^2-\tau(C).
\end{equation}
When $d=2d'$ is even, then one has
\begin{equation}
\label{Chern2}
\E_C=E_C(d'-1) \text{ and } c_1(\E_C)=-1, \ c_2(\E_C)=3(d')^2-3d'+1-\tau(C).
\end{equation}
Recall that the vector bundle $E_C$ is {\it stable} if and only if $\E_C$ has no sections, see \cite[Lemma 1.2.5]{OSS}, and in our case this is equivalent to $r=mdr(f)\geq d/2$,  see
also \cite[Proposition 2.4]{Se}.  Note that for $d$ even, $E_C$ is {\it semistable} if and only if it is stable, while for $d=2d'+1$ odd, $E_C$ is semistable if and only if $r \geq d'$. To see this, use the characterization of semistable rank 2 vector bundles on $\PP^n$ given by \cite[Lemma 1.2.5]{OSS}. Theorem \ref{thmNU} (2) and the formulas \eqref{Chern1} and \eqref{Chern2} imply that,
for a stable bundle $E_C$, one has
\begin{equation}
\label{Chern3}
 c_2(\E_C)=\nu(C).
\end{equation}
The following result was established in \cite{AD}, see Theorem 1.1, Proposition 3.1 and Proposition 3.2.
\begin{thm} \label{thm1} 
 With the above notation,  set $r=mdr(f)$. Then the following hold, where the line $L_0$ is generic and the line $L$ is arbitrary.
 
\begin{enumerate} 

\item   $d_1^L +d_2^L=d-1$;

\item $d_1^{L_0} \geq d_1^{L}$;

 \item $\max(r-\nu(C),0) \leq d_1^{L} \leq  d_1^{L_0} =  \min(r,\lfloor (d-1)/2\rfloor)$ and $$0\leq o(L)\leq \min(r,\nu(C)).$$
 
\item  $(d-1)^2-d_1^{L_0} d_2^{L_0}=\tau(C) + \nu(C).$

 \end{enumerate}

\end{thm}

\begin{rk}
\label{rktop}
The above formulas for Chern classes imply that, for two plane curves $C:f=0$ and $C':f'=0$
with $\deg C= \deg C'$ and $\tau(C)=\tau(C')$, the associated bundles $E_C$ and $E_{C'}$ are topologically equivalent, see for instance \cite[Section 6.1]{OSS}. In particular, this applies to  $C,C'$ the pair of line arrangements constructed by Ziegler in \cite{Zi}, such that 
$\deg C= \deg C'=9$ and $\tau(C)=\tau(C')=42$, see \cite[Remark 8.4]{DHA}.
Since $5=mdr(f) < mdr(f')=6$, it follows that $E_C$ and $E_{C'}$ are non-isomorphic stable vector bundles, even though $C$ and $C'$ have the same combinatorics.
However, $E_C$ and $E_{C'}$ have the same generic splitting type $(d_1^{L_0}, d_2^{L_0})$,
as follows from Theorem \ref{thm1} (3) above.
It seems that {\it no similar example exists involving unstable vector bundles}.
\end{rk}

Let $\alpha_L$ be the defining equation of the line $L$ in $X$. Then one has an exact sequence
$$0 \to \OO_X(-1) \stackrel{\cdot \alpha_L}{\to} \OO_X \to \OO_L \to 0,$$
where the first non-trivial morphism is induced by multiplication by the linear form $\alpha_L$. Let $k$ be an integer and tensor the above exact sequence by the vector bundle $E_C(k)$. We get
$$0 \to E_C(k-1) \stackrel{\cdot \alpha_L}{\to} E_C(k) \to E_C(k)|_L \to 0,$$
with $E_C(k)|_L \simeq \OO_L(k-d_1^L) \oplus \OO_L(k-d_2^L)$, since we assume as in the Introduction that $E_C|_L \simeq \OO_L(-d_1^L) \oplus \OO_L(-d_2^L)$. 
Then we have the following.

\begin{prop}
The long exact sequence of cohomology groups of the short exact sequence above starts as follows:
\begin{eqnarray} \label{equa2} 
0 &\to& AR(f)_{k-1} \stackrel{\cdot \alpha_L}{\to} AR(f)_{k} \stackrel{\pi_L}{\to} H^0(L,\OO_L(k-d_1^L) \oplus \OO_L(k-d_2^L)) \\
&\to& N(f)_{k+d-2} \stackrel{\cdot \alpha_L}{\to}
N(f)_{k+d-1} \to \cdots. \nonumber
\end{eqnarray}
Moreover, for $k=-1$, the corresponding morphism $N(f)_{d-3} \stackrel{\cdot \alpha_L}{\to}
N(f)_{d-2}$ is injective and hence $d_1^L\geq 0$  for any line $L$. 
\label{key1}
\end{prop}

\proof
This is exactly as in the proof of \cite[Theorem 5.7]{DS14}.  The key point is the identification 
\begin{equation}
\label{H1}
H^1(X, E_C(k))=N(f)_{k+d-1}, 
\end{equation}
valid for any integer $k$,  for which we refer 
to \cite[Proposition 2.1]{Se}. For the last claim, note that $N(f)_{d-3} \subset S_{d-3}$ and
$N(f)_{d-2} \subset S_{d-2}$, as the Jacobian ideal is generated in degree $d-1$.
\endproof
Finally, recall the following result, saying that the Jacobian module $N(f)$ enjoys a weak Lefschetz type property, see \cite{DPop} for this result and \cite{H+1,H+2,IG} for Lefschetz properties of Artinian algebras. 

\begin{thm}
\label{thm2}
If $L_0: \al_{L_0}=0$ is a generic line in $X$, then the morphism
$$N(f)_{s-1} \stackrel{\cdot \alpha_{L_0}}{\to}
N(f)_{s},$$
induced by the multiplication by $\alpha_{L_0}$, is injective for $s < \lceil T/2 \rceil$, and surjective for $s \geq \lceil T/2 \rceil$.
\end{thm}

See  Corollary \ref{corLef} and the end of Remark \ref{rkB} for partial strong Lefschetz property of the Jacobian module $N(f)$, the second one coming from a result by K. Hulek.

\section{On the Hilbert function of the Jacobian module $N(f)$}

The study of the dimensions $h^1(X, E_C(k))$ or, equivalently, in view of \eqref{H1}, the study of the Hilbert function 
$$n(f)_k=\dim N(f)_k$$
of the Jacobian module $N(f)$, is a central question in the study of rank 2 (stable) vector bundles on $X$, see for instance \cite{Ha78,Ha80}.
One has the following result for the vector bundle $E_C$, in the stable situation, saying that, in the middle range, the points $(j,n(f)_j)$ lie on an {\it upward pointing parabola}.

\begin{thm}
\label{thmB2} If $r=mdr(f) \geq d/2$, then the following hold for
$$2d-4-r \leq j \leq d-2+r.$$
\begin{enumerate} 

\item  For $d=2d'+1$ odd,  one has $T=3(d-2)=6d'-3$ and
$$n(f)_j=3(d')^2-(j-3d'+2)(j-3d'+1)-\tau(C)=\nu(C)-(j-\lfloor T/2 \rfloor )(j-\lceil T/2 \rceil).$$

\item  For $d=2d'$ even, one has $T=3(d-2)=6d'-6$ and
$$n(f)_j=3(d')^2-3d'+1-(j-3d'+3)^2-\tau(C)=\nu(C)-\left(j-\frac{T}{2}\right)^2.$$

 \end{enumerate}
\end{thm}

\proof
The equality of the two formulas for $n(f)_j$ in both cases follows from the formulas for $\nu(C)$ given in Theorem \ref{thmNU} (2).
One can derive a proof for the first equality in (1) above using \cite[Theorem 7.4 (a)]{Ha78}, case $-t-2\leq l\leq t-1$, and for the first equality in (2) using \cite[Theorem 7.4 (b)]{Ha78}, case $-t-1 \leq l \leq t-1$.
We present below an alternative proof.
 First we check both formulas for a smooth curve $C_F:f_F=0$, where $f_F=x^d+y^d+z^d$, when 
 $N(f)=M(f)$, and hence $n(f)_k=m(f)_k$ for all $k$.
The formulas for these dimensions are given for instance in \cite[Proposition 2.1]{St}, see in particular the explicit form for $n=2$ given just after the proof.
Consider now the general case, and note that
$$n(f)_j= m(f)_j-\dim (S/\hat J_f)_j.$$
By the definition of the coincidence threshold $ct(f)$, one has
$m(f)_j=m(f_F)_j$ for all $j \leq ct(f)$,
and 
$$ct(f)=d-2+mdr'(f) \geq d-2+r,$$
where $mdr'(f)$ is the minimal degree of a syzygy in $AR(f)$
which is not in the submodule $KR((f) \subset AR(f)$ generated by the Koszul relations $(f_y,-f_x,0)$, $(f_z,0,-f_x)$ and $(0,f_z,-f_y)$, see
\cite[Formula (1.3)]{DBull}.
On the other hand, one has
$$\dim (S/\hat J_f)_j=\tau(C)$$
for $j \geq T-ct(f)=3(d-2)-(d-2+mdr'(f))=2d-4-mdr'(f)$, see \cite[Proposition 2]{DBull}, and hence
in particular for $j \geq 2d-4-r$. This completes the second proof of Theorem \ref{thmB2}.

\endproof

The case of unstable rank 2 vector bundle on $X$ does not seem to have been considered until now. In this case, and assuming $C$ is not free, we have the following result, saying that, in the middle range, the points $(j,n(f)_j)$ lie on a {\it horizontal line segment, with a one-unit drop at the extremities}.

\begin{thm}
\label{thmB1} If $r=mdr(f)<d/2$  and $e$ is an integer such that $0 \leq e \leq 2$, then the following holds
$$n(f)_{d+r+e-5}= n(f)_{2d-r-e-1}=\nu(C)-\frac{(e-2)(e-3)}{2}+\alpha(C,e) ,$$
where $\alpha(C,e)\geq 0$.
In particular, we have the following.
\begin{enumerate} 

\item  For $e=2$, we get $\alpha(C,e)=0$ and
$$n(f)_j=\nu(C)\text{ for any } j \in [d+r-3,2d-r-3].$$

\item  For $e=1$, either $\alpha(C,e)=1$, and then
$C$ is free and  
$$n(f)_{d+r-4}=n(f)_{2d-r-2}=\nu(C)=0,$$
or $\alpha(C,e)=0$ and then
$C$ is not free and  $n(f)_{d+r-4}=n(f)_{2d-r-2}=\nu(C)-1$.

 \end{enumerate}
\end{thm}

\proof
Using the formula \eqref{Chern} above and the formulas in  \cite{Dmax}, (2.3), we get
\begin{equation} \label{ne1}
ar(f)_{k+1}+ar(f)_{d-5-k}+{d+k+2 \choose 2}-3{k+3 \choose 2}=n(f)_{d+k}+\tau(C),
\end{equation}
for any integer $k$. We set
$k+1=d-r-e$, for an integer $e \geq 0$, and we note that
\begin{equation} \label{ne2}
ar(f)_{d-r-e}=\dim S_{d-2r-e}\rho_1 +\alpha(C,e)={d-2r-e+2 \choose 2}+\alpha(C,e),
\end{equation}
for some integer $\alpha(C,e)\geq 0$, if $r=mdr(f)$ and we assume that $2r \leq d$, $e \leq 2$.
Since $d-5-k=r+e-4 \leq r-2$, we see that $ar(f)_{d-5-k}=0$ and a direct computation transforms the equation \eqref{ne1} into
\begin{equation} \label{ne3}
n(f)_{2d-r-e-1}+\tau(C)-\alpha(C,e)=(d-1)^2-r(d-r-1)-\frac{(e-2)(e-3)}{2}.
\end{equation}
Use next the  formula for $\nu(C)$ in Theorem \ref{thmNU} (1) and the well known duality result for $N(f)$ implying that
$n(f)_j=n(f)_{T-j}$ for any integer $j$, see \cite{Se}.
\endproof
The combination of Theorem \ref{thm2} and Theorem \ref{thmB1} yields the following
{\it partial strong Lefschetz property holds for the Jacobian module} $N(f)$.
\begin{cor} \label{corLef} 
If $r=mdr(f)<d/2$ and $L_0: \al_{L_0}=0$ is a generic line in $X$, then the morphism
$$N(f)_{p} \stackrel{\cdot \alpha_{L_0}^{q-p}}{\to}
N(f)_{q},$$
induced by the multiplication by $\alpha_{L_0}^{q-p}$, is an isomorphism for any
$$d+r-3 \leq p <q \leq 2d-r-3.$$

 \end{cor}

\section{Jumping lines and Lefschetz type properties for the Jacobian module}
The following result relates the splitting type of $E_C$ along a line $L: \al_L=0$, to the Lefschetz properties of the Jacobian module $N(f)$ with respect to the multiplication by $\al_L$.
\begin{prop} \label{propthmA} 
 For any line $L: \al_L=0$ in $X$, we have
 $d_1^L=\min \{mdr(f), k(f,L)\},$
 where $$k(f,L)= \min \{ k \in \N \ : \ N(f)_{k+d-2} \stackrel{\cdot \alpha_L}{\to}
N(f)_{k+d-1} \text{ is not injective } \}.$$
 \end{prop}

\proof If $k<\min \{mdr(f), k(f,L)\}$, the exact sequence \eqref{equa2} implies that $k<d_1^L$.
Hence $\min \{mdr(f), k(f,L)\} \leq d_1^L$.
If $k=mdr(f)$ or if $k=k(f,L)$, the same exact sequence implies $k \geq d_1^L$. Hence
$d_1^L \leq \min \{mdr(f), k(f,L)\}$, which proves our claim.
\endproof
The above proof also implies the following.

\begin{cor} \label{corA} 
Let $C:f=0$ be a reduced plane curve of degree $d$ and set $r=mdr(f)$. Then the following hold.
\begin{enumerate} 

\item  If $d_1^L=r$, then $L$ is not a jumping line,  $ar(f)_r \leq 2$, and the equality is possible only when $C$ is free with exponents $(d_1,d_1)$, and $d=2d_1+1$ is odd.

\item If $d_1^L <r< d_2^L$, then $ar(f)_r \leq r-d_1^L+1$.

\item If $d_2^L \leq r$, then $ar(f)_r \leq 2r-d+3$. 
 \end{enumerate}

  \end{cor}
  
The equality $ar(f)_r = 2r-d+3$  occurs when $C$ is a nearly free curve with $d=2d_1$ even and exponents $(d_1,d_1)$, and in many other cases, see Examples \ref{ex1} and \ref{ex2} below.
  \proof
  
  For the first claim note that $d_1^{L_0} \leq r$ by Theorem \ref{thm1} (4) or by Proposition \ref{propthmA}, and hence $L$ is not a jumping line. The inequaltiy $ar(f)_r \leq 2$ follows from the exact sequence \eqref{equa2} since $AR(f)_{r-1}=0$. If equality $ar(f)_r = 2$
  holds, it follows that $f$ has two linearly independent Jacobian syzygies, both of degree $r$.
  Hence the sum of their degrees is $2r=2d_1^L\leq d_1^L+d_2^L=d-1$. This is possible only when there are equalities everywhere and the curve $C$ is free with exponents $(r,r)$ by 
\cite[Lemma (1.1)]{ST}.  The remaining claims follow along the same lines.
  \endproof

\begin{rk} \label{rkA} 
If the morphism
$N(f)_{s-1} \stackrel{\cdot \alpha_L}{\to}
N(f)_{s}$
is not injective and $s \leq \lceil T/2 \rceil$, then the morphism
$N(f)_{s} \stackrel{\cdot \alpha_L}{\to}
N(f)_{s+1}$
is also not injective. Indeed, let $u \in N(f)_{s-1} $ be a non-zero element, such that $u\cdot \alpha_L=0$. Then, for a generic line $L_0= \al_{L_0}$, the element $u_0=u\cdot \al_{L_0} \in N(f)_{s}$ is non-zero by Theorem \ref{thm2}. On the other hand, it is clear that
$$u_0\cdot \alpha_L=  u\cdot \al_{L_0}\cdot \al _L= u\cdot \al_{L}\cdot \al _{L_0}=0.$$
In other words, the injective morphism $N(f)_{s-1} \stackrel{\cdot \alpha_{L_0}}{\to}
N(f)_{s}$ sends $K(\al_L)_{s-1}$ into $K(\al_L)_{s}$, where
$$K(\al_L)_{m}=\ker \{N(f)_{m} \stackrel{\cdot \alpha_L}{\to}
N(f)_{m+1}\}.$$
 \end{rk}
 
Now we investigate the jumping lines of $E_C$, namely the lines $L$ in $X$ such that $d_1^L<d_1^{L_0}$. Any line $L$ in $X$ corresponds clearly to a point in $\PP(S_1)$, corresponding to a defining linear form $\al_L$.
For any integer $k < mdr(f)$, consider the linear map
\begin{equation} \label{e2} 
\lambda _k:S_1 \to \Hom(N(f)_{d-2+k}, N(f)_{d-1+k}),
\end{equation}
sending a linear form $\al_L \in S_1$ to the morphism of multiplication by $\al_L$.
We assume that $d-2+k <T/2$, i.e. $k<(d-2)/2$, and hence 
$$n(f)_{d-2+k} \leq n(f)_{d-1+k},$$
by Theorem \ref{thm2}, where $n(f)_m=\dim  N(f)_{m}$ for any integer $m$.
Let
\begin{equation} \label{e3} 
\Sigma_k \subset  \Hom(N(f)_{d-2+k}, N(f)_{d-1+k}),
\end{equation}
denote the affine variety of linear maps which are not of maximal rank. 
Recall that
\begin{equation} \label{e3.5} 
\codim \Sigma_k =n(f)_{d-1+k}- n(f)_{d-2+k}+1,
\end{equation}
when $n(f)_{d-2+k}>0$, and $\Sigma_k=\emptyset$ when $n(f)_{d-2+k}=0$.

We define the $k$-th jumping locus of the curve $C:f=0$ to be the set
\begin{equation} \label{e4} 
V_k(C)= \{L \in \PP(S_1) \ : \ d_1^L \leq k\}.
\end{equation}

\begin{thm}
\label{thmB} If $k \geq mdr(f)$, then $V_k(C)=\PP(S_1)$. On the other hand,
for  $k <mdr(f)$, the following hold.
\begin{enumerate} 

\item  If $n(f)_{d-2+k}=0$, then $V_k(C)=\emptyset$.

\item  If $n(f)_{d-2+k}>0$, then $V_k(C)=(\lambda_k^{-1}(\Sigma_k) \setminus \{0\})/\C^*$ is a determinantal subvariety in
$\PP(S_1)=(S_1\setminus \{0\})/\C^*$.

\item $\emptyset=V_{-1}(C) \subset V_{0}(C) \subset ... \subset V_{d_1^{L_0}-1}(C) \subset V_{d_1^{L_0}}(C)= \PP(S_1)=\PP^2$.

\item If $\delta_k=n(f)_{d-1+k}=n(f)_{d-2+k}>0$, then $V_{k}(C)$ is a curve of degree at most $\delta_k$. 

\item If $\delta_k=n(f)_{d-1+k}=n(f)_{d-2+k}+1 >1$, then  $V_{k}(C)$ is either 1-dimensional, or 0-dimensional and $|V_{k}(C)|\le \delta_k(\delta_k-1)/2$ in this latter case. 

 \end{enumerate}
\end{thm}

\proof

The first two claims follow from the exact sequence \eqref{equa2}. The third claim follows from the inequality $d_1^{L}\leq d_1^{L_0}$, see Theorem \ref{thm1}, (2).
To prove (4), note that in this case $\Sigma_k$ is a hypersurface of degree $\delta_k$ and $0 \in \Sigma_k$.
Note that $\Lambda_k=\im \lambda_k$ is a linear space not contained in $ \Sigma_k$ by Theorem \ref{thm2}. It follows that $\lambda_k^{-1}(\Sigma_k)$ is a (possibly non-reduced) surface in $S_1=\C^3$, defined by a homogeneous polynomial of degree $\delta_k$.
The proof of the last claim is similar, in this case $\Sigma_k$ has codimension 2, and hence
$\lambda_k^{-1}(\Sigma_k)$ has codimension either 1 or 2, i.e. it cannot consist only of the origin $0$. When $\lambda_k^{-1}(\Sigma_k)$ has codimension 1, it consists of a number of lines, bounded by the degree of the determinantal variety $\Sigma_k$. This degree is known to be $\delta_k(\delta_k-1)/2$, see \cite[Example 19.10]{Har}.
\endproof

\begin{cor}
\label{corB} Let $C:f=0$ be a reduced plane curve of degree $d$ which is neither free nor nearly free, and assume that $r=mdr(f)$ satisfies $r<d/2$. Then the vector bundle $E_C$ is not stable, and it is semistable exactly when $d=2d'+1$ is odd and $r=d'$. Moreover, the following hold.
\begin{enumerate} 

\item  $d_1^{L_0}=r$ and hence $V_{r}(C)= \PP(S_1)=\PP^2$.

\item The set of jumping lines $V_{r-1}(C)$ is a curve of degree at most $\nu(C)$ in $\PP(S_1)$.

\item The set of jumping lines of order at least two $V_{r-2}(C)$ is  either 1-dimensional, or 0-dimensional and in this latter case $|V_{r-2}(C)|\le \nu(C)(\nu(C)-1)/2$ . 

 \end{enumerate}

\end{cor}

In Example  \ref{exS1} we have $d=5>4=2r$, $\nu(C)=3$ and $V_0(C)$ consists of 3 points,
hence the bound in Corollary \ref{corB} (3) is sharp in this case. The curve $V_{r-1}(C)$ is in fact a line arrangement in this case, as shown in Theorem \ref{thmD2} below.

\proof 

The first claim  in  Corollary \ref{corB} follows from Theorem \ref{thm1} (3), the second claim from 
Theorem \ref{thmB1} (1) and  Theorem \ref{thmB} (4) for $k=r-1$, and the final claim from
Theorem \ref{thmB1} (2) and  Theorem \ref{thmB} (5) for $k=r-2$.
\endproof

\begin{cor}
\label{corB1} Let $C:f=0$ be a reduced plane curve of degree $d$ which is not nearly free, and assume that $r=mdr(f)$ satisfies $r\geq d/2$. Then the vector bundle $E_C$ is  stable and the following hold.
\begin{enumerate} 

\item For $d=2d'+1$, one has $d_1^{L_0}=d'$ and set of jumping lines
$V_{d'-1}(C)$ is a curve of degree at most $\nu(C)$ in $\PP(S_1)$.

\item For $d=2d'$, one has $d_1^{L_0}=d'-1$ and set of jumping lines
$V_{d'-2}(C)$  is  either 1-dimensional, or 0-dimensional and in this latter case $|V_{d'-2}(C)|\le \nu(C)(\nu(C)-1)/2$ . 

 \end{enumerate}

\end{cor}
In Example \ref{exS}, for the Fermat quartic, we have $d=4<6=2r$ and set of jumping lines $V_0(C)$ is the union of 3 lines, hence a pure 1-dimensional variety.
In Example \ref{ex0}, we have $d=6=2r$ and set of jumping lines $V_1(C)$ is the union of a line and a point, hence it is 1-dimensional, but not pure 1-dimensional. On the other hand, in Example \ref{ex2}, we have $d=6<8=2r$, $\nu(C)=7$ and set of jumping lines $V_1(C)$ consists of 11 points.

\proof The first claim follows  from 
Theorem \ref{thmB2} (1) and  Theorem \ref{thmB} (4) for $k=d'-1$, and the final claim from
Theorem \ref{thmB2} (2) and  Theorem \ref{thmB} (5) for $k=d'-2$.
\endproof

 \begin{rk} \label{rkB} The claims above saying that some jumping sets $V_k(C)$ are pure 1-dimensional are related to Barth's Theorem (applied to our setting), see \cite{Ba},  \cite[Theorem 2.2.4]{OSS} as well as \cite{OSS}, pp. 118-119, saying that if $\E$ is a rank 2 vector bundle on $\PP^2$, which is semistable and has an even Chern class $c_1(\E)$, then the set of jumping lines of $\E$ is pure 1-dimensional.
In this situation, the equation of the curve $V_{d'-1}$ is given by the determinant of the mapping $N(f)_{3d'-2} \stackrel{\cdot \alpha_L}{\to}
N(f)_{3d'-1}$.

For semistable rank 2 vector bundles $\E$ with odd Chern class, i.e. $d=2d'$, the corresponding result 
to Barth's Theorem fails. An example of this situation for our bundles $E_C$ is given below in Example \ref{ex0}. In this case, K. Hulek has introduced in \cite{Hu} the notion of a {\it jumping line $L$ of the second kind}, which means that the mapping $N(f)_{3d'-4} \stackrel{\cdot \alpha^2_L}{\to}
N(f)_{3d'-2}$ has not maximal rank. Theorem 2.7 (2) implies that the set of jumping lines of the second kind is defined by the vanishing of the determinant $\Delta(a,b,c)$ of this latter mapping, 
regarded as a polynomial in the coefficients $a,b,c$ of the linear form $\al_L$. Hence the corresponding  jumping set $V_k(C)$ is a (possibly non-reduced) curve $C(E_C)$ of degree $2(\nu(C)-1)$, since this determinant $\Delta(a,b,c)$ is not identically zero by \cite[Theorem 3.2.2]{Hu}. See Example \ref{exS} below, the case of the Fermat quartic, for a situation where this curve (considered with reduced structure) has degree $<2(\nu(C)-1)$.

Note also that the non-vanishing of $\Delta(a,b,c)$ for a generic line $L$ implies that $$N(f)_{3d'-4} \stackrel{\cdot \alpha^2_L}{\to}
N(f)_{3d'-2}$$
is an isomorphism in this case, i.e. a {\it partial strong Lefschetz property holds for the Jacobian module} $N(f)$.

 \end{rk} 

\begin{ex} \label{exB} Let $C:f=0$ be a nearly free curve of degree $d$, with exponents $d_1\leq d_2$. When $d_1=d_2$ it is known that there are no jumping lines and the generic splitting type is $(d_1^{L_0},d_2^{L_0})=(d_1-1,d_1)$. 
The corresponding vector bundles $E_C$ is isomorphic to $T_X(-d_1-1)$, the shifted tangent bundle of $X$, see for details \cite{AD, MaVa}.
Consider now the case $d_1<d_2$, when it is known that  the generic splitting type is $(d_1^{L_0},d_2^{L_0})=(d_1,d_2-1)$, and a jumping line $L$ has a splitting type  $(d_1^{L},d_2^{L})=(d_1-1,d_2)$, see \cite{AD, MaVa}. 
Apply now Theorem \ref{thmB} to this situation. By \cite[Corollary 2.17]{DStRIMS}, we know that $n(f)_m=1$ if $d+d_1-3 \leq m \leq d+d_2-3$ and $n(f)_m=0$ otherwise.
If we apply Theorem \ref{thmB} (1) for $k=d_1-2< mdr(f)=d_1$, we have $V_k(C)=\emptyset$, i.e. the only possible splitting type is indeed $(d_1^{L},d_2^{L})=(d_1-1,d_2)$.
Apply now Theorem \ref{thmB} (4) for $k=d_1-1< mdr(f)=d_1$, and conclude that
$V_{d_1-1}(C)$ is a line, since $\delta_k=1$. A geometric description of this line was given in \cite{MaVa}, and a generalization of this result is discussed in our next section, see Theorem \ref{thmD}.
 \end{ex}

\section{Jumping lines and  the Bourbaki ideal of the syzygy module} 

Let $C:f=0$ be a reduced plane curve of degree $d$. For any choice of a nonzero syzygy
$\rho_1=(a_1,b_1,c_1) \in AR(f)_r$, where $r=mdr(f)$, we get a morphisms of graded $S$-modules
\begin{equation} \label{B1}
S(-r)  \xrightarrow{u} AR(f), \  u(h)= h \cdot \rho_1.
\end{equation}
For any syzygy $\rho=(a,b,c) \in AR(f)_m$, consider the determinant $\Delta(\rho)=\det M(\rho)$ of the $3 \times 3$ matrix $M(\rho)$ which has as first row $x,y,z$, as second row $a_1,b_1,c_1$ and as third row $a,b,c$. Then it turns out that $\Delta(\rho)$ is divisible by $f$, see \cite{Dmax}, and we define thus a new morphisms of graded $S$-modules
\begin{equation} \label{B2}
 AR(f)  \xrightarrow{v}  S(r-d+1)   , \  v(\rho)= \Delta(\rho)/f,
\end{equation}
and a homogeneous ideal $B(C,\rho_1) \subset S$ such that $\im v=B(C,\rho_1)(r-d+1)$.
The following result, except the claim (2), was stated for line arrangements in \cite[Proposition 2.1]{DStSuper}.

\begin{thm}
\label{thmC}
Let $C:f=0$  be a reduced plane curve of degree $d$ and
set $r=mdr(f)$. Then, for any choice of a nonzero syzygy
$\rho_1 \in AR(f)_r$, there is  an exact sequence
$$0 \to S(-r)  \xrightarrow{u} AR(f) \xrightarrow{v}  B(C,\rho_1)(r-d+1) \to  0,$$
and the following hold.
\begin{enumerate}

\item The ideal $B(C,\rho_1)$ is saturated, defines a subscheme $Z(C,\rho_1)=V(B(C,\rho_1))$ of $\PP^2$ of dimension at most 0, and its degree is given by
$$\deg B(C,\rho_1)=(d-1)^2-r(d-r-1)-\tau(C)=\tau_{max}(d,r)-\tau(C).$$

\item The ideal $B(C,\rho_1)$  and the codimension 2  subscheme $Z(C,\rho_1)$ are locally complete intersections.

\item The ideal $B(C,\rho_1)$  and the subscheme $Z(C,\rho_1)$ do not depend on the choice of $\rho_1$ when $\dim AR(f)_r=1$.

\item The curve $C$ is free if and only if $B(C,\rho_1)=S$.

\item The curve $C$  is nearly free if and only if the  subscheme $Z(C,\rho_1)$ is a reduced point $P(C, \rho_1)$ in $\PP^2$. The exact sequence and the point $P(C,\rho_1)$ are independent of $\rho_1$ when $2r<d$, i.e. when the exponents of the nearly free curve $C$ satisfy $r=d_1<d_2=d-r$.

\end{enumerate}
\end{thm} 

It follows that $B(C,\rho_1)$ is a Bourbaki ideal for the syzygy module $AR(f)$, see \cite{Bou}, Chapitre 7, \S 4, Thm. 6, as well as section 3 in \cite{SUV}. A similar construction for surfaces in $\PP^3$ was given in \cite{D0}. The dependence of the ideal $B(C,\rho_1)$ and of the scheme $Z(C,\rho_1)$ of the choice of the syzygy $\rho_1$ is illustrated in Example \ref{ex2}.
For $2r<d$, it follows from Theorem \ref{thmNU} (1) that 
$$\deg B(C,\rho_1)=\deg Z(C,\rho_1)=\nu(C).$$

\proof

We let the reader check that the proof given for  \cite[Proposition 2.1]{DStSuper} works as well in this more general setting. As for the new claim (2), we proceed as follows. If we sheafify the exact sequence of graded $S$-modules from Theorem \ref{thmC}, we get
 an exact sequence 
  \begin{equation} \label{ESVB0} 
0 \to \OO_X(-r)  \xrightarrow{\tilde u} E_C  \xrightarrow{ \tilde v} \I(r-d+1) \to 0.
 \end{equation}
 Here $\I$ is the sheaf ideal in $\OO_X$ associated to the Bourbaki ideal $B(C,\rho_1)$, and hence the support of $\OO_X / \I$ coincides with the support of the scheme $Z(C,\rho_1)$.
If $p$ belongs to this support, the surjectivity of $\tilde v_p$ implies that the corresponding  ideal $\I_p$ is generated by at most two elements. Indeed,
$E_{C,p}$ is a free $\OO_{X,p}$-module of rank 2. Since the scheme $Z(C,\rho_1)$ is 0-dimensional, this yields the claim (2).

\endproof

\begin{rk} \label{rkC} 
Note that the syzygy $\rho_1$ determines a section of the bundle $E_C(r)$, whose scheme of zeroes is exactly $Z(C,\rho_1)$. In particular, one has
$$\deg B(C,\rho_1)=\deg Z(C,\rho_1)=c_2(E_C(r)).$$
However, the {\it explicit construction} of the ideal $B(C,\rho_1)$ given above is useful, since it provides a simple method to obtain a minimal set of generators for this ideal $B(C,\rho_1)$.
Let $I(\rho_1)$ be the ideal in $S$ generated by the components $a_1,b_1,c_1$ of the syzygy $\rho_1$ and let $Z(I(\rho_1))$ be the 
corresponding subscheme in $\PP^2$. Then it is easy to see that the support $|Z(I(\rho_1))|$ of $Z(I(\rho_1))$ coincides with 
 the support $|Z(C,\rho_1)|$ of $Z(C,\rho_1)$ outside $C$. The example $C: f=x^5y^2z^2+x^9+y^9=0$, where $\rho_1=(-2xy^2z,0,9x^4+5y^2z^2)$, $|Z(I(\rho_1))|=\{(0:1:0), (0:0:1)\}$ and
$|Z(C,\rho_1)|=\{(0:1:0)\}$, shows that these two supports do not coincide in general. Note that in this example $r=mdr(f)=4$ and $ar(f)_4=1$, so the choice of $\rho_1$ is unique (up to a nonzero factor).

 \end{rk}

 \subsection{On lines avoiding the support of the jumping subscheme.} 
 
We discuss first the lines disjoint from the support of the jumping subscheme $Z(C,\rho_1)$.

\begin{thm} \label{thmD} 
Let $C:f=0$  be a reduced plane curve of degree $d$, set $r=mdr(f)$  and consider the subscheme  $Z(C,\rho_1)$ introduced above. Any line $L$ in $\PP^2$ which avoids the support of $Z(C,\rho_1)$ is not a jumping line if 
$2r \leq d.$
 More precisely, the (unordered) splitting type of $E_C$ along $L$ is $(r,d-1-r)$.
 \end{thm}
 
 \proof
 
If we tensor the exact sequence  \eqref{ESVB0} by $\OO_L$, for $L$ a line disjoint from the support of
$Z(C,\rho_1)$, we get the following exact sequence
 \begin{equation} \label{ESVB} 
0 \to \OO_L(-r)  \xrightarrow{\al } E_C|_L  \xrightarrow{\be} \OO_L(r-d+1) \to 0.
 \end{equation}
The isomorphism classes of such extensions of $\OO_L(r-d+1)$ by $\OO_L(-r)$ are classified by 
$$\Ext^1(\OO_L(r-d+1), \OO_L(-r))=\Ext^1(\OO_L, \OO_L(d-1-2r))=H^1(L,\OO_L(d-1-2r))=0,$$
see \cite[Section III.6]{Ha}, which proves our claim.
 \endproof
 
  \begin{cor} \label{corD} 
Let $C:f=0$  be a reduced plane curve of degree $d$, such that  $r=mdr(f) \leq d/2$. Then the set of jumping lines for the vector bundle $E_C$ is contained in a union of at most
$(d-1)^2-r(d-r-1)-\tau(C)$ lines in $\PP(S_1)$.
 \end{cor}

 \begin{rk} \label{rkD2} The condition $2r \leq d$ in Theorem \ref{thmD} is necessary, as Example \ref{ex1} below shows.  \end{rk} 
 
 \begin{thm} \label{thmBD} 
Let $C:f=0$  be a reduced plane curve of degree $d$ and consider the subscheme  $Z(C,\rho_1)$ introduced above. Then, if $r=mdr(f)>d/2$, the splitting type $(d_1^L,d_2^L)$ along any line $L$ in $\PP^2$ which avoids the support of $Z(C,\rho_1)$ satisfies 
$ d_1^ L\geq d-1-r$. In particular, if $2r-d \in \{1,2\}$, then $d_1^ L \in \{d_1^ {L_0}-1, d_1^ {L_0}\}$.
  \end{thm}
 Examples in the next section shows that this lower bound is sharp in many cases, e.g. in the situation of the last claim, both values for $d_1^L$ are obtained,
 see the final parts of Examples \ref{ex1} and \ref{ex2}.
 
 \proof
 We use the same notation as in the proof of Theorem \ref{thmD}. In the exact sequence \eqref{ESVB} we have $E_C|_L =\OO_L(-d_1^L)\oplus 
 \OO_L(-d_2^L)$. The surjective morphism $\be$ is induced by a pair of homogeneous polynomials $(A_1,A_2) \in S_{a_1}\times S_{a_2}$,
 where $a_i=r-d+1+d_i^L$ for $i=1,2$, satisfying the condition $G.C.D.(A_1,A_2)=1$. Indeed, at the level of sections, the morphism $\be$ is given by
 $$(s_1,s_2) \mapsto A_1s_1+A_2s_2.$$
 Note that $a_1 \leq a_2$. If $A_1 \ne 0$, then $a_1 \geq 0$, and this yields the claim of our Theorem. If $A_1=0$, it follows that $A_2$ is a non-zero constant, and hence $a_2=0$. This implies 
 $$d_2^L=d-1-r <\frac{d-1}{2},$$
 which is a contradiction. Indeed, $d_2^L \geq d_1^L$ implies 
 $$d_2^L\geq \frac{d-1}{2}.$$
 The last claim follows by checking that, in these two situations, one has
 $$d-r-1=d_1^ {L_0}-1.$$
   \endproof

 \subsection{On lines meeting the support of the jumping subscheme $Z(C,\rho_1)$} 
 
 Let $L$ be a line in $\PP^2$ such that $L \cap | Z(C,\rho_1)|=\{p_1,...,p_s\}$. For each such point $p_k$ we define its multiplicity as follows. Consider a system of local coordinates $(u,v)$ centered at $p_k$ such the equation of the line $L$ is given by $u=0$. The localized ideal $\I_{p_k} \subset \OO_{X,p_k}=\C\{u,v\}$, being a complete intersection, is generated by two analytic germs, say $g(u,v)$ and $h(u,v)$. Then we set
 $$ m_k =\dim _{\C} \frac{\C\{u,v\}}{(u,g(u,v),h(u,v))}=\dim _{\C} \frac{\C\{v\}}{(g(0,v),h(0,v))}.$$
Then clearly $1 \leq m_k<+ \infty$ and one has
$$\frac{\C\{u,v\}}{(u)} \otimes_{\C\{u,v\}} \frac{\C\{u,v\}}{(g(u,v),h(u,v))}=\frac{\C\{u,v\}}{(u,g(u,v),h(u,v))},$$
and hence the latter ring can be regarded as the local ring of the point $p_k$ in the scheme theoretic intersection $Z(C,\rho_1) \cap L$.
The ideal $\I_{p_k}=(g(u,v),h(u,v)) \subset \OO_{X,p_k}$, being a complete intersection, we have a free resolution
$$ 0 \to  \OO_{X,p_k} \to  \OO_{X,p_k}^2 \to \I_{p_k} \to 0,$$
where the non-trivial morphisms are given by the pair $(g(u,v),h(u,v))$.
When we  tensor by $\OO_{L,p_k}$, we get the following exact sequence
$$  \OO_{L,p_k} \to  \OO_{L,p_k}^2 \to \I_{p_k} \otimes \OO_{L,p_k} \to 0,$$
and the corresponding morphisms are given by the pair $(g(0,v),h(0,v)) \ne (0,0)$. It follows that the first morphism is injective, and up-to a change of basis in $\OO_{L,p_k}^2=\C\{v\}^2$ is given by the pair
$(v^{m_k},0)$. It follows that 
$$\I_{p_k} \otimes \OO_{L,p_k}=\C\{v\} \oplus \frac{\C\{v\}}{(v^{m_k})}.$$
 If we tensor now the exact sequence \eqref{ESVB0} by $\OO_L$, we get, keeping track of the twists and using the above local computations, the following result. When the points $p_k \in Z(C,\rho_1) \cap L$ are all simple points, then this result is already in \cite{FV1}, see equation (7).
 
 \begin{prop} \label{propE} 
 With the above notation, there is an exact sequence
 $$0 \to \OO_L(-r) \to E_C|_L \to \OO_L(r-d+1-m_L) \oplus \left(\oplus_{k=1,s}\frac{\OO_{L,p_k}}{M_{p_k}^{m_k}}\right) \to 0,$$
where $m_L=\sum_{k=1,s}m_k$ and  $M_{p_k} \subset \OO_{L,p_k}$ denotes the corresponding maximal ideal.
 \end{prop}
 
Using this Proposition and its notation, we can prove the following result.
 
\begin{thm} \label{thmD2} 
Let $C:f=0$  be a reduced plane curve of degree $d$, set $r=mdr(f)$  and consider the subscheme  $Z(C,\rho_1)$ introduced above. Any line $L$ in $\PP^2$ which meets the support of $Z(C,\rho_1)$ is  a jumping line if 
$2r \leq d-1.$
 More precisely, the splitting type of $E_C$ along $L$ is $(r-m_L,d-1-r+m_L)$ or, equivalently, the order of the jumping line $L$ is given by $o(L)=m_L \leq r$. Moreover, the set of jumping lines $V_{r-1}(C)$ is a line arrangement consisting of at most $\nu(C)$ lines, dual to the support of the subscheme  $Z(C,\rho_1)$.
 \end{thm} 
 \proof
 
 It is clear that the splitting type of $E_C$ along $L$ is $(r-h,d-1-r+h)$ for some $0 \leq h \leq r$. If $0 \leq h <m_L$, then we have $-r+h \geq r-d+1-h >r-d+1-m_L$, and hence there is no surjective morphism from 
$E_C|_L$  to $\OO_L(r-d+1-m_L)$, which is a contradiction in view of Proposition \ref{propE}. It follows that $h \geq m_L$. Assume now that $h>m_L$. Then $-r>r-d+1-h$, and hence the first nontrivial morphism in the exact sequence from Proposition \ref{propE} is given by a pair $(H,0)$, where $H$ is a homogeneous polynomial of degree $h>m_L$.
This implies that the torsion part of the cokernel of this morphism has dimension equal to $h>m_L$, a contradiction.

Since the degree of the subscheme $Z(C,\rho_1)$ is $\nu(C)$ for $2r\leq d-1$ by Theorem \ref{thmC} (1) and Theorem \ref{thm1} (3) and (5),  the last claim follows as well.
\endproof

A computation of the splitting type using this approach can be seen in Example \ref{exS1}.
 
 \begin{rk} \label{rkD1} 
 (i) Note that the unstable rank 2vector bundles on $X=\PP^2$ have been studied by Schwarzenberger in \cite{Sch1} under the name of {\it almost decomposable} vector bundles. The equivalence of the two notions follows for instance from \cite[Theorems 1.2.9 and 1.2.10]{OSS}.
Schwarzenberger  has shown that for such a vector bundle, the set of jumping lines is a union of {\it pencils}, that is lines in the dual projective plane, see \cite[Proposition 10]{Sch1}. Since $E_C$ is unstable exactly when $2r<d$, our Theorem \ref{thmD2} can be regarded as a refinement of Schwarzenberger's result for the bundles $E_C$.

\noindent (ii) The example of a nearly free curve $C$ with exponents $(d_1,d_1)$ discussed in Example \ref{exB}, when there are no jumping lines but the scheme $Z(C,\rho_1)$ consists of a simple point, shows that a line $L$ meeting the support of $Z(C,\rho_1)$ may not be a jumping line if $r=mdr(f) \geq d/2$.
A similar situation is described in Examples  \ref{ex1} and \ref{ex0} below. Note that Example \ref{ex0} shows that  the set of jumping lines described in Corollary \ref{corD} is not necessarily pure 1-dimensional, i.e. it may consists of lines and isolated points, when $r= d/2$, unlike the case $r<d/2$ covered by Theorem \ref{thmD2}.
 \end{rk}

 \section{Some examples} 
 
 First we consider the smooth curves.
 \begin{ex} \label{exS}
Let $C:f=0$ be a smooth curve of degree $d \geq 3$. Then $r=mdr(f)=d-1$ and the graded $S$-module $AR(f)$ is generated by the Koszul type syzygies
$$\rho_1=(f_y,-f_x,0), \ \rho_2=(f_z,0,-f_x) \text{  and  } \rho_3=(0,f_z,-f_y).$$
With this choice, the Bourbaki ideal $B(C,\rho_1)$ is spanned by
$v(\rho_2)=d \cdot f_x$ and $v(\rho_3)=d \cdot f_y$, hence it is a global complete intersection. For the Fermat type curve 
$$C:f_F=x^d+y^d+z^d=0,$$
the support of the scheme $Z(C,\rho_1)$ is the multiple point $p=(0:0:1)$. The line $L:z=0$ does not pass through this point and Proposition \ref{propthmA} implies that $d_1^L=0=d-r-1$, i.e. for this line we get equality in the inequality given by Theorem \ref{thmBD}.

\medskip
 {\bf Case $d=2d'+1$ odd.} Then Corollary \ref{corB1} implies that the set of jumping lines
$V_{d'-1}(C)$ is a curve in $\PP(S_1)$. The geometry of these curves $V_{d'-1}(C)$ depends on the equation $f$.
For instance, in the case of a plane cubic
$$C: f=x^3+y^3+z^3+3t xyz=0, \text{ where  } t \in \C,  \  t^3\ne -1,$$
an easy direct computation shows that
\begin{equation} \label{sym1.5} 
V_{d'-1}(C): t(a^3+b^3+c^3)+(2-t^3)abc=0,
\end{equation} 
where $(a:b:c)$ are the coordinates on $\PP(S_1)$. In other words, the jumping variety $V_{d'-1}(C)$ determines the complex structure of $C$ up to finite indeterminacy in this case. This is related to Barth's result in 
\cite[Theorem 2 and Section 7]{Ba}, saying that a rank 2 stable vector bundle on $\PP^2$ with even second Chern class is determined by the associated net of quadrics, having the curve $V_{d'-1}(C)$ as its discriminant.

\medskip
 {\bf Case $d=2d'$ even.} Then Corollary \ref{corB1} implies that the set of jumping lines
 $V_{d'-2}(C)$ is nonempty. For $f=x^4+y^4+z^4$ and using the usual monomial bases for $N(f)=M(f)$, we get $V_0(C):abc=0$, hence the union of 3 lines. In particular, $V_0(C)$ is pure 1-dimensional in this case. Note that the determinant of the maping $N(f)_{2} \stackrel{\cdot \alpha^2_L}{\to}
N(f)_{4}$, where $\alpha_L=ax+by+cz$, is given by $a^4b^4c^4$. Hence the curve of jumping lines of second order $C(E_C)$ is given by the equation $a^4b^4c^4=0$, and hence its support
coincides with $V_0(C)$ in this case. In other words, we have equality in \cite[Proposition 9.1]{Hu}.

 \end{ex}

 The computations in the following examples  were all  done using the computer algebra software SINGULAR, see \cite{Sing}. The Chern classes of $E_C$ can be computed in each case using \eqref{Chern} above, since we give in each example the corresponding global Tjurina number $\tau(C)$.
 
  \begin{ex} \label{exS1}
 Let $C:f=0$, where $f=x^5+y^5+(x^4+y^4)z$. Then $d=5$, $\tau(C)=9$, and 
 $r=mdr(f)=2$, therefore the bundle $E_C$ is semistable. Theorem \ref{thm1} (3) implies that the corresponding generic splitting type of $E_C$ is $(d_1^{L_0},d_2^{L_0})=(2,2)$.
 The Jacobian ideal $J_f$ is spanned by $f_x,f_y,f_z$, and its saturation  $\wJ_f$ is spanned by $ x^3, y^3$. The only non-zero dimensions $n(f)_m$ are in this case
 $n(f)_4=n(f)_5=3$ and $n(f)_3=n(f)_6=2$. Moreover, a vector space basis of $N(f)_3$ (resp. of $N(f)_4$) is given by
 $x^3, y^3$ (resp. $x^4, x^3y,xy^3$). With respect to these bases, the multiplication
 $\{N(f)_{3} \stackrel{\cdot \alpha_L}{\to}
N(f)_{4}\}$, where $ \al_L=ax+by+cz$, is given by
$$(ax+by+cz)\cdot x^3=(a-\frac{5c}{4})x^4+bx^3y$$
and
$$(ax+by+cz)\cdot y^3=(\frac{5c}{4}-b)x^4+axy^3.$$
It follows  that $V_0(C)$ consists of  3 points, namely
$(0:0:1), (0:5:4),(5:0:4)$. Since $\nu(C)=3$, it follows that
 we have equality in Corollary \ref{corB} (3), hence the bound there is sharp in this situation. Similarly, a basis for $N(f)_5$ is given by $x^5,x^3y^2,x^2y^3$ and the multiplication
 $\{N(f)_{4} \stackrel{\cdot \alpha_L}{\to}
N(f)_{5}\}$ is given by
$(ax+by+cz)\cdot x^4=(a+b-\frac{5c}{4})x^5,$
$(ax+by+cz)\cdot x^3y=(a-\frac{5c}{4}-b)x^5+bx^3y^2$ and 
$(ax+by+cz)\cdot xy^3=-(b-\frac{5c}{4})x^5+ax^2y^3.$
It follows  that $V_1(C)$ consists of  3 lines, namely $\LL_1:a=0$, $\LL_2:b=0$ and $\LL_3:4(a+b)-5c=0$.
The $S$-module $AR(f)$ has 4 generating syzygies, of degrees 2,4,4,4 and a direct computation shows that the scheme $Z(C,\rho_1)$, which does not depend on the choice of the syzygy $\rho_1$, consists of the simple points $P_1=(1:0:0)$, $P_2=(0:1:0)$ and $P_3=(4:4:-5)$. It follows that the line $\LL_j \subset \PP(S_1)$ above consists of all the lines in $\PP^2$ passing through the point $P_j$, for $j=1,2,3$.
Note that the corresponding lines $L=L_{i,j}$ in $V_0(C)$ pass through the points $P_i,P_j$ in the support of $Z(C,\rho_1)=\{P_1,P_2,P_3\}$, and one has $m_L=r=2$ in this case, as predicted by Theorem \ref{thmD2}. More precisely, one has $L_{1,2}:z=0$, $L_{1,3}:5y+4z=0$ and $L_{2,3}:5x+4z=0$.

 \end{ex}

 \begin{ex} \label{ex1}
 Let $C:f=0$, where $f=2x^5+2y^5+5x^2y^2z$. Then $d=5$, $\tau(C)=10$, and we see that the $S$-module $AR(f)$ is generated by 4 syzygies $\rho_i$, $i=1,...,4$, all of  degree 
 $r=mdr(f)=3$. Hence Theorem \ref{thm1} (3) implies that the corresponding generic splitting type of $E_C$ is $(d_1^{L_0},d_2^{L_0})=(2,2)$.
 The Jacobian ideal $J_f$ is spanned by $f_x,f_y,f_z$, and its saturation  $\wJ_f$ is spanned by $f_x,f_y,f_z, x^3y, xy^3$. The only non-zero dimensions $n(f)_m$ are in this case
 $n(f)_4=n(f)_5=2$. Moreover, a vector space basis of $N(f)_4$ (resp. of $N(f)_5$) is given by
 $x^3y, xy^3$ (resp. $x^4y, xy^4$). With respect to these bases, the multiplication
 $\{N(f)_{4} \stackrel{\cdot \alpha_L}{\to}
N(f)_{5}\}$, where $ \al_L=ax+by+cz$, is given by
$$(ax+by+cz)\cdot x^3y=ax^4y-cxy^4$$
and
$$(ax+by+cz)\cdot xy^3=-cx^4y+bxy^4.$$
Using Theorem \ref{thmB} (4) for $k=1$, we get that $V_1(C)$, the set of jumping lines for $E_C$, is the smooth conic $Q:ab-c^2=0$ in $\PP(S_1)$.

Hence in this case we have
$$\emptyset=V_{-1}(C) =V_{0}(C) \subset  V_{1}(C)=Q \subset V_{2}(C)= \PP(S_1).$$
Indeed, Theorem \ref{thmB} (1) implies that $V_{0}(C)=\emptyset$.
If we choose
$$\rho_1=(0,x^2y,-2(y^3+x^2z))\in AR(f)_3,$$
then the corresponding Bourbaki
$B(C,\rho_1)$ is the ideal $(xz,y^2,xy)$, and hence the scheme $Z(C,\rho_1)$ consists of two points, a simple one at $(1:0:0)$, given in local coordinates by an ideal $(u,v)$, and a double point at $(0:0:1)$, given in local coordinates by an ideal $(u,v^2)$.

Among the lines on $Q$, only the lines $x=0$ and  $y=0$ meet the support of $Z(C,\rho_1)$.
For the other lines in $Q$, the bound given by Theorem \ref{thmBD} is $d_1^L \geq d-r-1=1$.
In fact, we have equality, hence this bound is sharp in this situation. 
 \end{ex} 

\begin{rk} \label{rkE} 
 
 The smooth conic $Q$ above  is one of the smooth degree $n$ curves occurring as jumping loci, predicted by Barth for stable rank 2 vector bundles $\E$ on $\PP^2$, with $c_1(\E)=0$ and $c_2(\E)=n$, see \cite[ Application 1, section 5.4]{Ba}. Indeed, note that the normalization of our vector bundle $E_C$ is $\E_C=E_C(2)$ and it satisfies $c_1(\E_C)=0$ and $c_2(\E_C)=2$.  Similar remarks apply for the cubic curve in \eqref{sym1.5},
 which is smooth for $t^3 \notin \{-1,0, 8\}$.

\end{rk} 

  \begin{ex} \label{ex0}
   Let $C:f=0$, where $f=(x^2+y^2)^3+(y^3+z^3)^2$, i.e. $C$ is a Zariski sextic with 6 cusps on a conic. Then $d=6$, $\tau(C)=12$, and  we see that the $S$-module $AR(f)$ is generated by 4 syzygies $\rho_i$, $i=1,...,4$, of degrees 
 $r=mdr(f)=3=d_1<d_2=d_3=d_4=5$. Hence Theorem \ref{thm1} (4) implies that the corresponding generic splitting type of $E_C$ is $(d_1^{L_0},d_2^{L_0})=(2,3)$.
 The Jacobian ideal $J_f$ is spanned by $f_x,f_y,f_z$, and its saturation  $\wJ_f$ is spanned by $g=y^3+z3$ and $h=(x^2+y^2)^2$. The only non-zero dimensions $n(f)_m$ are in this case $n(f)_3=n(f)_9=1$,
 $n(f)_4=n(f)_8=4$, $n(f)_5=n(f)_7=6$ and $n(f)_6=7$. Moreover, a vector spaces basis of $N(f)_5$ (resp. of $N(f)_6$) is given by
 $$ x^2g,y^2g, xyg,xzg,yzg,zh,$$
 and respectively by
 $$ x^3g,x^2yg, y^3g,x^2zg,y^2zg,xyzg,z^2h.$$
  With respect to these bases, the multiplication
 $\{N(f)_{5} \stackrel{\cdot \alpha_L}{\to}
N(f)_{6}\}$, where $\al_L=ax+by+cz$, is given by the matrix
\begin{center}
$$M(L)=\left(
  \begin{array}{ccccccc}
     a & 0& 0 & 0 & 0 & 0 \\
     b & 0& a& 0 & 0 & 0 \\
     0& b & 0 & 0 & 0& 0 \\
   c& 0 & 0 & a & 0 & 0 \\
    0 & c& 0 & 0 & b & -b \\
   0 & 0 & c & b& a& 0  \\
   0 & 0 & 0 & 0 & 0& c  \\
    \end{array}
\right)$$
\end{center}
Using Theorem \ref{thmB} (5) for $k=1$, we get that $V_1(C)$, the set of jumping lines for $E_C$, is the set of lines $L$ such that $\rank M(L) < 6$. A direct computation, shows that $V_1(C)$ consists of the line $\LL:a=0$ and one points, namely 
$P_1=(1:0:0)$.
 A vector basis for $N(f)_4$ is given by $xg,yg,zg,h$, and using the given bases, the multiplication
 $\{N(f)_{4} \stackrel{\cdot \alpha_L}{\to}
N(f)_{5}\}$, where $\al_L=ax+by+cz$, is given by the matrix
\begin{center}
$$M'(L)=\left(
  \begin{array}{ccccccc}
     a & 0& 0  & 0\\
     0 & b& 0& -b\\
     b& a & 0 & 0 \\
   c& 0 & a & 0 \\
    0 & c& b & 0\\
   0 & 0 & 0 & c \\
  
    \end{array}
\right)$$
\end{center}
Using Proposition \ref{propthmA}, it follows that $V_{0}(C)$ is the set of lines $L$ such that
$\rank M'(L)<4$, which implies that $V_0(C)=\{P_1,P_2,P_3\}$, where $P_1$ is as above, $P_2=(0:1:0)$ and $P_3=(0:0:1)$.
Hence in this case we have
$$\emptyset=V_{-1}(C) \subset V_{0}(C)=\{P_1,P_2,P_3\} \subset  V_{1}(C)=\{P_1\} \cup \LL \subset V_{2}(C)= \PP(S_1).$$
Since $ar(f)_3=1$, there is essentially a unique choice
$$\rho_1=(yz^2,-xz^2, xy^2)\in AR(f)_3.$$
Then the corresponding Bourbaki
$B(C,\rho_1)$ is the ideal $(xy^2,xz^2, yz^2)$, and hence the scheme $Z(C,\rho_1)$ consists of three points, one nonreduced at $p_1=(1:0:0)$, given in local coordinates by an ideal $(u^2,v^2)$, the other nonreduced at $p_2=(0:1:0)$, given in local coordinates by an ideal $(u,v^2)$, and a reduced point at $p_3=(0:0:1)$ given by $(u,v)$.  Note that the line $\LL$ consists of all the lines passing through the point $p_1$, the line $L_1:x=0$, corresponding to the point $P_1$, passes through the points $p_2$ and $p_3$,  the line $L_2:y=0$,
 corresponding to the point $P_2$, passes through the points $p_1$ and $p_3$, and
the line $L_3:z=0$,
 corresponding to the point $P_3$, passes through the points $p_1$ and $p_2$.  None of the points $p_i$ is situated on the sextic $C$.
  \end{ex}

  \begin{ex} \label{ex2}
 Let $C:f=0$, where $f=x^6+y^6+3x^2y^2z^2$. Then $d=6$, $\tau(C)= 12$,  and  we see that the $S$-module $AR(f)$ is generated by 5 syzygies $\rho'_i$, $i=1,...,5$,  of degrees 
 $r=mdr(f)=4=d_1=d_2 <d_3=d_4=d_5=5$, see their expressions given below. Hence Theorem \ref{thm1} (4) implies that the corresponding generic splitting type of $E_C$ is $(d_1^{L_0},d_2^{L_0})=(2,3)$.
 The Jacobian ideal $J_f$ is spanned by $f_x,f_y,f_z$, and its saturation  $\wJ_f$ is spanned by $f_x,f_y,f_z, x^3y, x^2y^2,xy^3$. The only non-zero dimensions $n(f)_m$ are in this case
 $n(f)_4=n(f)_8=3$, $n(f)_5=n(f)_7=6$ and $n(f)_6=7$. Moreover, a vector spaces basis of $N(f)_5$ (resp. of $N(f)_6$) is given by
 $$ xy^4,x^2y^3,x^3y^2,x^4y, xy^3z,x^3yz,$$
 and respectively by
 $$ xy^5,x^2y^4,x^3y^3, x^4y^2,x^5y, xy^4z,x^4yz.$$
  With respect to these bases, the multiplication
 $\{N(f)_{5} \stackrel{\cdot \alpha_L}{\to}
N(f)_{6}\}$, where $\al_L=ax+by+cz$, is given by the matrix
\begin{center}
$$M(L)=\left(
  \begin{array}{ccccccc}
     b & 0& 0 & 0 & 0 & -c \\
     a & b& 0& 0 & 0 & 0 \\
     0& a & b & 0 & 0& 0 \\
   0& 0 & a & b & 0 & 0 \\
    0 & 0& 0 & a & -c & 0 \\
   c & 0 & 0 & 0& b& 0  \\
   0 & 0 & 0 & c & 0& a  \\
    \end{array}
\right)$$
\end{center}
Using Theorem \ref{thmB} (5) for $k=1$, we get that $V_1(C)$, the set of jumping lines for $E_C$, is the set of lines $L$ such that $\rank M(L) < 6$. A direct computation shows that $V_1(C)$ consists of the following 11 points in $\PP(S_1)$:
$$P_1=(1:1:1), \ P_2=(1:1:-1), \ P_3=(1:-1:1), \ P_4=(-1:1:1),$$
$$P_5=(1:0:0), \ P_6=(0:1:0), \ P_7=(0:0:1), \ P_8=(\al^2:\al:1),$$
$$P_9=(\al:\al^2:1), \ P_{10}=(\be^2:\be:1), \ P_{11}=(\be:\be^2:1),$$
where $\al^2+\al+1=0$ and $\be^2-\be+1=0$. A vector basis for $N(f)_4$ is given by $xy^3,x^2y^2,x^3y$, and using the given bases, the multiplication
 $\{N(f)_{4} \stackrel{\cdot \alpha_L}{\to}
N(f)_{5}\}$, where $\al_L=ax+by+cz$, is given by the matrix
\begin{center}
$$M'(L)=\left(
  \begin{array}{ccccccc}
     b & 0& 0  \\
     a & b& 0\\
     0& a & b  \\
   0& 0 & a  \\
    c & 0& 0 \\
   0 & 0 & c  \\
  
    \end{array}
\right)$$
\end{center}
Using Proposition \ref{propthmA}, it follows that $V_{0}(C)$ is the set of lines $L$ such that
$\rank M'(L)<3$, which implies that $V_0(C)=P_7=(0:0:1)$.
Hence in this case we have
$$\emptyset=V_{-1}(C) \subset V_{0}(C)=\{P_7\} \subset  V_{1}(C)=\{P_j \ : \ j=1,...,11\} \subset V_{2}(C)= \PP(S_1).$$
The software SINGULAR gives the following minimal system of generators for the graded $S$-module $AR(f)$:
$$\rho'_1=(0,-x^2yz, y^4+x^2z^2),  \  \rho'_2=(-xy^2z,0, x^4+y^2z^2), \ \rho'_3=(xyz^3,-x^4z, x^2y^3-yz^4),$$
$$\rho'_4=(-y^4z,xyz^3,x^3y^2-xz^4) \text{ and } \rho'_5=(-y^5-x^2yz^2,x^5+xy^2z^2,0).$$
Since now $ar(f)_4=2$, there are several choices for the syzygy $\rho_1$ in Theorem \ref{thmC}. We discuss three choices.
\medskip

{\bf Choice 1.} If we choose $\rho_1=\rho'_1$, 
then the corresponding Bourbaki ideal
$B(C,\rho_1)$ is spanned by $g_2=v(\rho'_2)=-xyz$, $g_3=v(\rho'_3)=xz^3$, $g_4= v(\rho'_4)=-y^3z$ and $g_5=v(\rho'_5)=-y^4-x^2z^2$, where $v$ is the morphism defined in \eqref{B2}. Hence the scheme $Z(C,\rho_1)$ consists of two points, both nonreduced, one at $p_1=(1:0:0)$, given in local coordinates $u,v$ by an ideal $(uv, v^2+u^4)$, and another at $p_2=(0:0:1)$, given in local coordinates by an ideal $(u,v^3)$. 

\medskip

{\bf Choice 2.} If we choose $\rho_1=\rho'_2$, 
then the corresponding Bourbaki ideal
$B(C,\rho_1)$ is  spanned by $h_1=v(\rho'_1)=xyz$, $h_3=v(\rho'_3)=x^3z$, $h_4= v(\rho'_4)=-yz^3$ and $h_5=v(\rho'_5)=-x^4-y^2z^2$. Hence the support of the scheme $Z(C,\rho_1)$ consists of two points, one at $q_1=(0:1:0)$, and another at $p_2=(0:0:1)$, the same point as in Choice 1.

\medskip

{\bf Choice 3.} If we choose $\rho_1=\rho'_1+t\rho'_2$, where $t \in \C^*$, 
then the corresponding Bourbaki ideal 
$B(C,\rho_1)$ is  spanned by $k_1=v(\rho'_1)=txyz$, $k_2=v(\rho'_2)=-xyz$,
$k_3=v(\rho'_3)=xz(z^2+tx^2)$, $k_4= v(\rho'_4)=-yz(y^2+tz^2)$ and $k_5=v(\rho'_5)=
-y^4-x^2z^2-t(x^4+y^2z^2)=-y^2(y^2+tz^2)-x^2(tx^2+z^2)$. If we take $t=-s^4$ for $s \in \C^*$,
then the support of the scheme $Z(C,\rho_1)$ consists of the following  9 points:

\begin{itemize}
		\item[i)]  $z_j(s)=(\epsilon_j:s:0)$ for $j=1,2,3,4$, where $\epsilon_j$ are the four roots of  $\epsilon^4=1$;
\item[ii)] $z_j(s)=(0:s^2: (- 1)^j)$, where $j=5,6$;

\item[iii)] $z_j(s)=((- 1)^j:0:s^2)$, where $j=7,8$ and

\item[iv)] $z_9(s)=p_2=(0:0:1)$.
\end{itemize}
Theorem \ref{thmC} (1) tells us that $\deg B(C,\rho_1)=9$, and hence all these points $z_j(s)$ are simple points. When $s \to 0$, we see that the points 6 points $z_j(s)$ for $j\in \{1,2,3,4,7,8\}$
converge to the point $p_1$, and the 2 points $z_j(s)$ for $j\in \{5,6\}$
converge to the point $p_2=z_9(s)$. Similarly, when $|s| \to +\infty$, the 6 points
$z_j(s)$ for $j\in \{1,2,3,4,5,6\}$
converge to the point $q_1$, and the 2 points $z_j(s)$ for $j\in \{7,8\}$
converge to the point $p_2=z_9(s)$. Moreover, the line $L_7:z=0$, corresponding to the point $P_7$, contains the 4 points $z_j(s)$ for $j\in \{1,2,3,4\}$ for any $s$, the maximal number of collinear points among the points $z_j(s)$.
Note that the line $L_1:x+y+z=0$ is disjoint from the support of the scheme $Z(C,\rho_1)$ for most choices of $\rho_1$, and the bound given by Theorem \ref{thmBD} is $d_1^L \geq d-r-1=1.$ In fact, we have equality, hence this bound is sharp in this situation as well. 

 \end{ex} 
 
  \begin{rk} \label{rkFIN} 
(i)  In Example \ref{ex2}, the stable vector bundle $E_C$ admits a {\it unique jumping line $P_7$ of maximal order} $o(P_7)=2$. Note that condition (a) in \cite[Theorem 6.2]{Ha80} is not fulfilled, hence we cannot use this result of Hartshorne to deduce the unicity of a jumping line of maximal order.

\noindent (ii) A twist of the stable vector bundle $E_C$ in Example \ref{ex2} admits a section with 9 simple zeros $z_j(s)$ as explained in the third choice for $\rho_1$. However, the set of jumping lines does not coincides with the set of all lines passing through these points, and the line $P_7: z=0$ of maximal order 2 contains 4 of these points $z_j(s)$. This should be compared with \cite[Theorem 2.2.5]{OSS} and the discussion preceding it.

 \end{rk}


\begin{thebibliography}{00}

\bibitem{AD}  T. Abe, A. Dimca, On the splitting types of bundles of logarithmic vector fields along plane curves, Internat. J. Math. 29 (2018), no. 8, 1850055, 20 pp.

%\bibitem{A1}T. Abe, Chambers of 2-affine arrangements and freeness of 3-arrangements. 
%\textit{J. Alg. Combin}., \textbf{38} (2013), no. 1, 65--78. 

%\bibitem{A}T. Abe, Roots of characteristic polynomials and and intersection points of line arrangements. \textit{J. Singularities}, \textbf{8} (2014), 100--117.

%\bibitem{AN}T. Abe and Y. Numata, Exponents of $2$-multiarrangements and multiplicity lattices. \textit{J. Alg. Combin.}, \textbf{35} (2012),  no. 1, 1--17.


%\bibitem{ATW}T. Abe, H. Terao and M. Wakefield, The characteristic polynomial of a multiarrangement.\textit{Adv. in Math.}, \textbf{215} (2007), 825--838.

\bibitem{B+} E. Artal Bartolo, L. Gorrochategui, I. Luengo, A. Melle-Hern\' andez,
On some conjectures about free and nearly free divisors, in: {\it Singularities and Computer Algebra, Festschrift for Gert-Martin Greuel on the Occasion of his 70th Birthday}, pp. 1--19, Springer (2017)

\bibitem{Ba} W. Barth, Moduli of vector bundles on the projective plane, \textit{Invent. Math.} \textbf{42} (1977), 63--91.



\bibitem{Bou} N. Bourbaki, \textit{Alg\`ebre Commutative}, Chapitres I-IX, Hermann, 1961-1983.


\bibitem{CHMN} D. Cook, B. Harbourne, J. Migliore, U. Nagel,  Line arrangements and configurations of points with an unexpected geometric property. \textit{Compositio Math.} \textbf{154}(2018), 2150--2194.


%\bibitem{DS}G. Denham and M. Schulze,Complexes, duality and Chern classes of logarithmic forms along hyperplane arrangements.\textit{Arrangements of hyperplanes—Sapporo} 2009,  27-–57, Adv. Stud. Pure Math., 62, Math. Soc. Japan, Tokyo, 2012.

\bibitem
{Sing} { W. Decker, G.-M. Greuel, G. Pfister \and H. Sch{\"o}nemann.} \newblock {\sc Singular} {4-0-1} --- {A} computer algebra system for polynomial computations, available at {http://www.singular.uni-kl.de} 
(2014).

\bibitem{DBull}  A. Dimca, Syzygies of Jacobian ideals and defects of linear systems, \textit{Bull. Math. Soc. Sci. Math. Roumanie},  \textbf{56(104)} (2013),191--203.

\bibitem{D0}
A. Dimca, Freeness versus maximal degree of the singular subscheme for surfaces in $P^3$, 
\textit{Geom. Dedicata} \textbf{183} (2016), 101--112.

%\bibitem{D1} A. Dimca, Curve arrangements, pencils, and Jacobian syzygies,  \textit{Michigan Math. J.} \textbf{66} (2017), 347--365.

\bibitem{Dmax}
A. Dimca, 
Freeness versus maximal global Tjurina number for plane curves.
\textit{Math. Proc. Cambridge Phil. Soc.}, \textbf{163} (2017), 161--172.

\bibitem{Drcc} A. Dimca, On rational cuspidal plane curves, and the local cohomology of Jacobian rings, arXiv:1707.05258, to appear in \textit{Commentarii Mathematici Helvetici.}


\bibitem{DHA}  A. Dimca,   {\em Hyperplane Arrangements: An Introduction}, Universitext, Springer, 2017.

\bibitem{DPop} A. Dimca, D. Popescu, 
Hilbert series and Lefschetz properties of dimension one almost complete intersections, \textit{Comm. Algebra} \textbf{44} (2016), 4467--4482.


\bibitem{DS14} A. Dimca, E. Sernesi,  Syzygies and logarithmic vector fields along plane curves,
\textit{Journal de l'\'Ecole polytechnique-Math\'ematiques} \textbf{1} (2014), 247-267.


\bibitem{DStRIMS}
A. Dimca, G. Sticlaru, 
Free and nearly free curves vs. rational cuspidal plane curves, \textit{Publ. RIMS Kyoto Univ. }\textbf{54} (2018), 163--179.

\bibitem{DStSuper}
A. Dimca, G. Sticlaru, On supersolvable and nearly supersolvable line arrangements,
Journal of Algebraic Combinatorics,
https://doi.org/10.1007/s10801-018-0859-6.


%\bibitem{duPCTC} A.A. du Plessis,  C.T.C. Wall, Application of the theory of the
%discriminant to highly singular plane curves, \textit{Math. Proc. Camb.
%Phil. Soc.},  \textbf{126} (1999), 259-266. 






\bibitem{FV1}  D. Faenzi, J. Vall\`es, 
{Logarithmic bundles and line arrangements, an approach via the standard construction}, \textit{J. London.Math.Soc.}
\textbf{90} (2014), {675--694}.

\bibitem{H+1} T. Harima, J. Migliore, U. Nagel and J. Watanabe, { The weak and strong
Lefschetz properties for artinian K-algebras}, J. Algebra \textbf{262} (2003), 99--126.

\bibitem{H+2} T. Harima,  T. Maeno, H.  Morita, Y.  Numata, A.  Wachi, J.  Watanabe, {\it The Lefschetz properties}, 
Lecture Notes in Mathematics 2080,
Springer, Heidelberg, 2013.

\bibitem{Har} J. Harris,  \emph{Algebraic geometry: a first course},  Graduate Texts in Mathematics 133,  Springer-Verlag, New York, Heidelberg, Berlin, 1992.



\bibitem{Ha} R. Hartshorne,  \emph{ Algebraic geometry}, Graduate Texts in Mathematics 52, Springer-Verlag, New York, Heidelberg, Berlin, 1977.

\bibitem{Ha78} R. Hartshorne, Stable vector bundles of rank 2 on $\PP^3$, \textit{Math. Ann.} \textbf{238} (1978),  229--280.

\bibitem{Ha80} R. Hartshorne, Stable reflexive sheaves, \textit{Math. Ann.} \textbf{254} (1980), 121--176.

\bibitem{Hu}  K. Hulek, Stable rank 2 vector bundles on $\PP^2$ with $c_1$ odd, \textit{Math. Ann}. \textbf{242} (1979), 241--266.

\bibitem{IG} G. Ilardi, { Jacobian ideals, arrangements and Lefschetz properties}, \textit{J. Algebra}  \textbf{508} (2018), 418--430. 

\bibitem{MaVa} S. Marchesi, J. Vall\` es, Nearly free curves and arrangements: a vector bundle point of view, arXiv:1712.04867.


%\bibitem{MuS}  M.\ Musta\c t\u a, H.\ Schenck,
%The module of logarithmic p-forms of a locally free arrangement,
%\textit{J.\ Algebra} \textbf{241} (2001), 699--719.


\bibitem{OSS} C. Okonek, M. Schneider, H. Spindler: \emph{Vector Bundles on Complex Projective Spaces}. Progress in Math. n. 3, Birkhauser (1980).

%\bibitem{OT} P. Orlik, H. Terao :  \emph{Arrangements of Hyperplanes}. Grundlehren Math. Wiss., \textbf{300}, Springer-Verlag, Berlin  (1992)

\bibitem{S}
K. Saito, 
Theory of logarithmic differential forms and logarithmic vector fields.
\textit{J. Fac. Sci. Univ. Tokyo} \textbf{27} (1980), 265--291.  

\bibitem{Sch1} R.L.E. Schwarzenberger,  Vector bundles on algebraic surfaces, \textit{Proc. London Math. Soc.}  \textbf{11} (1961), 601--622. 


%\bibitem{Sch2} R.L.E. Schwarzenberger, Vector bundles on the projective plane. Proc. London Math. Soc. 11 (1961), 623--640.

\bibitem{Se} E. Sernesi:  The local cohomology of the jacobian ring, \textit{Documenta Mathematica},  \textbf{19} (2014), 541-565. 

\bibitem{Serre} J.-P. Serre, Sur les modules projectifs, S\'em. Dubreil-Pisot 1960/61, expos\'e 2.

\bibitem{ST} A. Simis, S.O. Toh\u aneanu, Homology of homogeneous divisors, \textit{Israel J. Math.} \textbf{200} (2014), 449-487.

\bibitem{SUV} A. Simis, B. Ulrich, W.V. Vasconcelos, Rees algebras of modules, \textit{Proc. LMS} \textbf{87} (2003), 610--646.

\bibitem{St} G. Sticlaru, Some criteria to check if a projective hypersurface is smooth or singular, British Journal of Mathematics \& Computer Science, \textbf{4} (2014), 924--932. 

%\bibitem{T1}
%H. Terao, 
%Arrangements of hyperplanes and their freeness I, II. 
%\textit{J. Fac. Sci. Univ. Tokyo} \textbf{27} (1980), 293--320.   

%\bibitem{T2}
%H. Terao, 
%Generalized exponents of a free arrangement of hyperplanes and
%Shephard-Todd-Brieskorn formula. \textit{Invent. Math}. 
%\textbf{63}  (1981),
%159--179.

%\bibitem{W}A. Wakamiko, On the exponents of 2-multiarrangements. 
%\textit{Tokyo J. Math}. \textbf{30} (2007), no. 1, 99--116. 

%\bibitem{WY}
%M. Wakefield and S. Yuzvinsky,
%Derivations of an effective divisor on the complex projective line.\textit{Trans. Amer. Math. Soc}. \textbf{359} (2007), 4389--4403. 

%\bibitem{Y1}M. Yoshinaga, Characterization of a free arrangement and conjecture of
%Edelman and Reiner. \textit{Invent. Math.} \textbf{157} (2004), no. 2,
%449--454.

%\bibitem{Y2} M. Yoshinaga,  On the freeness of 3-arrangements. 
%\textit{Bull. London Math. Soc.} \textbf{37} (2005), no. 1, 126--134. 



%\bibitem{Y3} M. Yoshinaga, Freeness of hyperplane arrangements and related topics. 
%\textit{Annales de la Facult\`{e} des Sciences de Toulouse } \textbf{23} (2014), no. 2, 483--512. 

%\bibitem{Z} G. M. Ziegler,  Multiarrangements of hyperplanes and their freeness.  Singularities (Iowa City, IA, 1986),  345--359, Contemp. Math., {\bf 90}, Amer. Math. Soc., Providence, RI, 1989. 

\bibitem{Zi} G. Ziegler, Combinatorial construction of logarithmic differential forms, Adv. Math. 76 (1989), 116-154.


\end{thebibliography}
\end{document}